\documentclass[a4paper, 12pt]{article}

\usepackage{amsfonts, amsmath, amsthm}
\usepackage{amssymb}
\usepackage{latexsym}
\usepackage[dvips]{graphicx}
\usepackage{color}

\theoremstyle{plain}
\newtheorem{thm}{Theorem}[section]
\newtheorem{prp}{Proposition}[section]
\newtheorem{lem}{Lemma}[section]

\theoremstyle{definition}

\theoremstyle{remark}
\newtheorem{rmk}{Remark}[section]

\numberwithin{equation}{section}

\newcommand{\R}{\mathbb{R}}
\newcommand{\C}{\mathbb{C}}

\newcommand{\cc}[1]{\overline{#1}}
\newcommand{\op}[1]{\mathcal{#1}}

\newcommand{\pa}{\partial}
\newcommand{\eps}{\varepsilon}
\newcommand{\jb}[1]{\langle #1 \rangle}

\DeclareMathOperator{\realpart}{\rm Re}
\DeclareMathOperator{\imagpart}{\rm Im}

\begin{document}
\title{
Large time asymptotics for a cubic nonlinear Schr\"odinger system 
in one space dimension\\
 }

\author{
          Chunhua Li \thanks{
              Department of Mathematics, College of Science, 
              Yanbian University. 
              977 Gongyuan Road, Yanji, 
              Jilin 133002, China. 
              (E-mail: {\tt sxlch@ybu.edu.cn})
             }
          \and
          Yoshinori Nishii\thanks{
              Department of Mathematics, Graduate School of Science, 
              Osaka University. 
              1-1 Machikaneyama-cho, Toyonaka, 
              Osaka 560-0043, Japan. 
              (E-mail: {\tt y-nishii@cr.math.sci.osaka-u.ac.jp})             }
           \and  
          Yuji Sagawa \thanks{
              Department of Mathematics, Graduate School of Science, 
              Osaka University. 
              1-1 Machikaneyama-cho, Toyonaka, 
              Osaka 560-0043, Japan. 
             }
           \and  
          Hideaki Sunagawa \thanks{
              Department of Mathematics, Graduate School of Science, 
              Osaka University. 
              1-1 Machikaneyama-cho, Toyonaka, 
              Osaka 560-0043, Japan. 
              (E-mail: {\tt sunagawa@math.sci.osaka-u.ac.jp})
             }
}

\date{\today }   
\maketitle

\noindent{\bf Abstract:}\ We consider a two-component system of cubic 
nonlinear Schr\"odinger equations in one space dimension. 
We show that each component of the solutions to this system behaves like a 
free solution in the large time, 
but there is a strong restriction between the profiles of them. 
This turns out to be a consequence of non-trivial long-range 
nonlinear interactions.
\\

\noindent{\bf Key Words:}\ 
Nonlinear Schr\"odinger system, asymptotically free, long-range interaction.
\\

\noindent{\bf 2010 Mathematics Subject Classification:}\ 
35Q55, 35B40.

\section{Introduction }  \label{sec_intro}

This paper deals with the global in time behavior of solutions 
$u=(u_1(t,x), u_2(t,x))$ to 
\begin{align}
\left\{\begin{array}{ll}
\begin{array}{l}
 \op{L} u_1=- i|u_2|^2 u_1,\\
 \op{L} u_2=- i|u_1|^2 u_2,
 \end{array}
 & (t,x )\in (0,\infty)\times \R,
 \end{array}\right.
\label{eq}
\end{align}
with the initial condition 
\begin{align}
u_j(0,x)=\varphi^0_j(x),
\qquad  x \in \R,\ j=1,2,
\label{data}
\end{align}
where $i=\sqrt{-1}$, $\op{L}=i\pa_t+(1/2)\pa_{x}^2$, and  
$\varphi^0=(\varphi^0_1(x), \varphi_2^0(x))$ is a given $\C^2$-valued 
function of $x\in \R$ which belongs to an appropriate weighted Sobolev space 
and satisfies a suitable smallness condition.

First of all, let us summarize the backgrounds briefly.
As is well-known, cubic nonlinearity gives a critical situation 
when we consider large time behavior of solutions to the nonlinear 
Schr\"odinger equation in $\R$. 
In general, cubic nonlinearity should be regarded as a long-range 
perturbation. 
For example, according to Hayashi--Naumkin \cite{HN}, the small data solution $u(t,x)$ to 
\begin{align}
 \op{L} u = \lambda |u|^2 u
 \label{nls_1}
\end{align}
with $\lambda \in \R\backslash\{0\}$ behaves like 
\[
 u(t,x)=\frac{1}{\sqrt{it}} \alpha^{\pm} (x/t) 
 e^{i\{\frac{x^2}{2t}  -\lambda  |\alpha^{\pm}(x/t)|^2 \log t \}}
 +o(t^{-1/2})
\quad \mbox{in} \ \ L^{\infty}(\R_x)
\]
as $t\to \pm \infty$, where  $\alpha^{\pm}$ is a suitable $\C$-valued function 
on $\R$. 
An important consequence of this asymptotic expression is that the solution to 
\eqref{nls_1} decays like $O(|t|^{-1/2})$ uniformly in $x\in \R$, 
while it does not behave like the free solution (unless $\lambda=0$). 
In other words, the additional logarithmic correction in the phase 
reflects a typical long-range character of the cubic nonlinear Schr\"odinger 
equations in one space dimension. 
If $\lambda\in \C$ in \eqref{nls_1}, another kind of long-range effect can be observed. For instance, according to \cite{Shi} 
(see also \cite{KitaShim}, \cite{JJL}, \cite{HLN}, etc.), 
the small data solution $u(t,x)$ to \eqref{nls_1} decays like 
$O(t^{-1/2}(\log t)^{-1/2})$ in $L^{\infty}(\R_x)$ as $t\to +\infty$ 
if $\imagpart\lambda<0$. 
This gain of additional logarithmic time decay should be interpreted as 
another kind of long-range effect. 
There are various extensions of these results. 
In the previous works \cite{LiS1} and \cite{LiS2}, 
several structural conditions on the nonlinearity have been introduced 
under which the small data global existence holds for a class of cubic 
nonlinear Schr\"odinger systems in $\R$, 
and large time asymptotic behavior of the global solutions have also been 
investigated (see also \cite{Kim}, \cite{SakSu} and the references cited 
therein for related works). 
We do not state these conditions here, but we only point out that 
the small data global existence for \eqref{eq} follows from the results of 
\cite{LiS1} and \cite{Kim} but the large time asymptotic behavior of 
solutions is not covered by these results. 
We note that the system \eqref{eq} possesses two conservation laws 
\[
 \frac{d}{dt}\Bigl( \|u_1(t)\|_{L^2}^2+ \|u_2(t)\|_{L^2}^2 \Bigr)=
-2\int_{\R} |u_1(t,x)|^2 |u_2(t,x)|^2\, dx
\]
and
\begin{align}\label{CL2}
 \frac{d}{dt}\Bigl( \|u_1(t)\|_{L^2}^2 - \|u_2(t)\|_{L^2}^2\Bigr)=0.
\end{align}
However, these are not enough to say something about the large time 
asymptotics for $u(t)$, and this is not trivial at all. 
To the authors' knowledge, 
there are no previous results which cover the asymptotic behavior of 
solutions to \eqref{eq}--\eqref{data}.

Our motivation of considering \eqref{eq} comes from the recent work~\cite{NS}, 
in which the system of semilinear wave equations
\begin{align}
\left\{\begin{array}{ll}
\begin{array}{l}
 (\pa_t^2-\Delta) v_1=- |\pa_t v_2|^2 \pa_t v_1,\\
 (\pa_t^2-\Delta) v_2=- |\pa_t v_1|^2 \pa_t v_2,
 \end{array}
 & (t,x )\in (0,\infty)\times \R^2,
 \end{array}\right.
\label{wave}
\end{align}
has been treated in connection with 
the Agemi-type structural condition (that is a kind of weak null conditions). 
From the viewpoint of conservation laws, there are a lot of similarities 
between \eqref{eq} and \eqref{wave}. It has been shown in \cite{NS} that 
global solutions to \eqref{wave} with small data behaves like 
solutions to the free wave equations, but there is a strong restriction 
in the profiles. Although the approach of \cite{NS} does not use 
the conservation laws directly, it may be natural to expect that 
an analogous phenomenon can be observed for solutions to  \eqref{eq}. 
The aim of the present paper is to reveal it. 

Before stating the main result, let us introduce some notations. 
For $s$, $s' \ge 0$, we denote by $H^s$ the $L^2$-based 
Sobolev space of order $s$, and the weighted Sobolev space $H^{s,s'}$ 
is defined by 
$\{\phi\in L^2\, |\, \jb{\,\cdot\, }^{s'} \phi \in H^{s}\}$ 
equipped with the norm 
$\|\phi\|_{H^{s,s'}}=\|\jb{\,\cdot\,}^{s'} \phi\|_{H^s}$, 
where $\jb{x}=\sqrt{1+x^2}$. The Fourier transform $\op{F}$ is defined by
\[
 \op{F}\phi(\xi) =\hat{\phi}(\xi)=\frac{1}{\sqrt{2\pi}} \int_{\R} 
 e^{-iy\xi} \phi(y)\, dy,\quad \xi\in \R.
\]
We also set $\op{U}(t)=\exp(i\frac{t}{2}\pa_x^2)$, so that 
$\op{U}(t)\phi=:w(t)$ solves the free Schr\"odinger equation 
$\op{L}w=0$ with $w(0)=\phi$.

The main result is as follows. 
\begin{thm}\label{thm_initial}
Suppose that $\varphi^0=(\varphi_1^0,\varphi_2^0)\in H^{2}\cap H^{1,1}$ 
and 
$\eps=\|\varphi^0\|_{H^2\cap H^{1,1}}$ is suitably small. 
Let $u=(u_1,u_2)\in C([0,\infty); H^2\cap H^{1,1})$ 
be the solution to \eqref{eq}--\eqref{data}. Then 
there exists $\varphi^+=(\varphi_1^+,\varphi_2^+)\in L^2$ with 
$\hat{\varphi}^+=(\hat{\varphi}_1^+,\hat{\varphi}_2^+)\in L^{\infty}$ 
such that 
\[
 \lim_{t\to +\infty} \|u_j(t)-\op{U}(t)\varphi_j^+\|_{L^2}=0,\quad j=1,2.
\]
Moreover we have 
\begin{align}\label{decoupling_ini}
\hat{\varphi}_1^+(\xi)\cdot \hat{\varphi}_2^+(\xi)=0,\quad \xi \in \R.
\end{align}
\end{thm}

\begin{rmk}\label{rem1} 
We emphasize that \eqref{decoupling_ini} should be regarded as a consequence 
of non-trivial long-range nonlinear interactions 
because such a phenomenon does not occur in the usual short-range situation. 
To complement this point, we will give auxiliary results on 
the final state problem for \eqref{eq} in Appendix \ref{sec_app_b}.\\
\end{rmk}

\begin{rmk} \label{rem2}
In the case where $u_1=u_2$, the system \eqref{eq} 
is reduced  to the single equation \eqref{nls_1} with $\lambda=-i,$ 
so  we can adapt the result of \cite{HLN} to see that 
$\|u(t)\|_{L^2}$ converges to $0$ as $t\to +\infty$ without restrictions on
the size of the initial data. 
However, this is an exceptional case. 
We are interested in general situations of \eqref{eq}--\eqref{data}. 
When $\|\varphi_1^0\|_{L^2}\ne \|\varphi_2^0\|_{L^2}$ in \eqref{data}, 
 it follows from the conservation law \eqref{CL2} and the $L^2$-unitarity of 
$\op{U}(t)$ and $\op{F}$ that at least one of $\hat{\varphi}_1^+$ or 
$\hat{\varphi}_2^+$ does not identically vanish. This implies that solutions 
$u=(u_1,u_2)$ to \eqref{eq}--\eqref{data} do not converge to $0$ as 
$t\to +\infty$ in $L^2$ for {\em generic} initial data being suitably small. 
In this sense, our problem is much more delicate than the single case 
\eqref{nls_1} with $\lambda=-i$.\\
\end{rmk}

\begin{rmk} \label{rem3}
It is worthwhile to note that the presence of $-i$ in the right-hand sides of \eqref{eq} 
is essential for our result. Indeed, if we drop $-i$ from the right-hand sides of \eqref{eq}
(that is, $\op{L} u_1= |u_2|^2 u_1$ and $\op{L} u_2= |u_1|^2 u_2$), we can show that 
the solutions have logarithmic phase corrections as in the single case 
\eqref{nls_1} with $\lambda \in \R\backslash \{0\}$ 
(see e.g. \cite{V} for detail).\\
\end{rmk}

\begin{rmk}\label{rem4}
The above theorem concerns only the forward Cauchy problem 
(i.e., for $t>0$). In the backward case (or, equivalently, 
in the case where $-i$ is replaced by $+i$ in the right-hand sides 
of \eqref{eq}), even the small data global 
existence may fail in general. See \cite{Sag} and the references cited 
therein for more information and the related works on this issue.\\
\end{rmk}

We close this section with the contents of this paper. 
We introduce some preliminary lemmas in 
Section \ref{sec_prelin}. Theorem  \ref{thm_initial} is proved in 
Section~\ref{sec_proof_initial}. 
In Appendix~\ref{sec_app_a}, we give a proof of a technical lemma. 
Appendix~\ref{sec_app_b} is devoted to auxiliary 
observations on the final state problem for \eqref{eq}.

\section{Preliminaries}  \label{sec_prelin}

In this section, we collect 
several identities and inequalities which are useful in the proof of 
Theorem~\ref{thm_initial}. In what follows, we will denote various 
positive constants by the same letter $C$, which may vary from one line to 
another. 
We set $\op{J}= x+it \pa_{x}$. It is well-known that 
$[\op{L},\op{J}]=0$ and $[\pa_x,\op{J}]=1$, where 
$[\cdot, \cdot]$ stands for the commutator, that is, 
$[\op{A},\op{B}]=\op{A}\op{B}-\op{B}\op{A}$ for two linear operators 
$\op{A}$ and $\op{B}$. Also we have 
\begin{gather}
\op{J}=\op{U}(t)x \op{U}(-t),\label{J_to_U}\\
 \|\phi\|_{L^{\infty}}
 \le 
 {C}{t^{-1/2}} \|\phi\|_{L^2}^{1/2} \|\op{J}\phi\|_{L^2}^{1/2}, 
\quad t>0,\label{kla_ineq}
\end{gather}
and 
\[
 \op{J}(\phi_1 \phi_2 \cc{\phi_3})
= 
(\op{J} \phi_1) \phi_2 \cc{\phi_3} 
+
\phi_1 (\op{J}\phi_2) \cc{\phi_3} 
-
\phi_1 \phi_2 \cc{\op{J}\phi_3}.
\]
Let $u=(u_1,u_2)$ be a smooth solution to \eqref{eq}--\eqref{data} on 
$[0, \infty) \times \R$. 
We define $\alpha=(\alpha_1,\alpha_2)$ by
\begin{align}
\alpha_j(t,\xi)=\op{F} \Bigl[\op{U}(-t)u_j(t,\cdot)\Bigr](\xi)
\label{def_alpha}
\end{align} 
for $j=1,2$. Then from (\ref{eq}) it follows that 
\begin{align}
\pa_{t} \alpha_1
=
-i\op{F}\op{U}(-t)\op{L}u_1
=
-\op{F}\op{U}(-t)(|u_2|^2u_1)
=
-\frac{1}{t} |\alpha_2|^2 \alpha_1   + R_1,
\label{profile_1}
\end{align}
where
\begin{align*}
R_1
=
\frac{1}{t} |\alpha_2|^2 \alpha_1
-\op{F}\op{U}(-t)\bigl[ |u_2|^2 u_1\bigr]. 
\end{align*}
Similarly we have 
\begin{align}\label{profile_2}
\pa_{t} \alpha_2=-\frac{1}{t} |\alpha_1|^2 \alpha_2 + R_2,
\end{align}
where
\begin{align*}
R_2
=
\frac{1}{t} |\alpha_1|^2 \alpha_2
-\op{F}\op{U}(-t)\bigl[ |u_1|^2 u_2\bigr].
\end{align*}
Concerning estimates for $R=(R_1,R_2)$, we have the following estimate.
\begin{lem}\label{remainder}
Let $R$ be as above. For $t\geq 1$, we have 
\begin{align*}
 |R(t,\xi)|
\le 
\frac{C}{t^{5/4}\jb{\xi}}
\bigl(\|u(t)\|_{H^1} + \|\op{J}u(t)\|_{H^1}\bigr)^{3}.
\end{align*}
\end{lem}
This estimate is not a new one (see e.g. Lemma~5.2 in \cite{LiS1}). 
For the convenience of the readers, we will give a proof in 
Appendix \ref{sec_app_a}.

Next we recall the  basic decay estimates for global solutions $u$ to 
\eqref{eq}--\eqref{data}.
From the argument of \cite{LiS1}, we already know the following result.
\begin{lem}\label{lem_est_u}
Let $0<\gamma< 1/12$. 
Suppose that $\eps=\|\varphi^0\|_{H^2\cap H^{1,1}}$ is suitably small. 
Then the solution $u$ to \eqref{eq}--\eqref{data} satisfies
\begin{align*}
\|u(t)\|_{H^2} +\|\op{J}u(t)\|_{H^1}\le C \eps \jb{t}^{\gamma}
\end{align*}
for $t\ge 0$ and
\begin{align}
|\alpha(t,\xi)| \le C \eps \jb{\xi}^{-1}
\label{apriori2}
\end{align}
for $t\ge 0$, $\xi \in \R$, where $\alpha$ is given by \eqref{def_alpha}.
\end{lem}

As a by-product of Lemmas~\ref{remainder} and \ref{lem_est_u}, we have
\begin{align}
|R(t,\xi)|
&\le 
\frac{C\eps^3}{t^{5/4-3\gamma} \jb{\xi}}
\label{R}
\end{align} 
for $t \geq 1$.
Roughly speaking, this means that the evolution of 
$\alpha=(\alpha_1,\alpha_2)$ may be governed by 
\[
\pa_t\alpha_1=-\frac{1}{t}|\alpha_2|^2 \alpha_1, 
\quad
\pa_t\alpha_2=-\frac{1}{t}|\alpha_1|^2 \alpha_2
\] 
up to the harmless remainders. 
We also note that $u(t)=\op{U}(t)\op{F}^{-1}\alpha(t)$. 
This point of view, whose original idea goes back to 
Hayashi--Naumkin~\cite{HN}, is the key of our approach.

\section{Proof of Theorem~\ref{thm_initial}}  \label{sec_proof_initial}

This section is devoted to the proof of Theorem~\ref{thm_initial}. 
The main step is to show the following. 

\begin{prp} \label{prp_asym_prof} 
Let $\alpha=(\alpha_1(t,\xi),\alpha_2(t,\xi))$ be given by \eqref{def_alpha} 
for the solution $u=(u_1,u_2)$ to \eqref{eq} satisfying the assumptions of 
Theorem~\ref{thm_initial}. There exists 
$\alpha^{+}=(\alpha_1^+(\xi), \alpha_2^+(\xi))\in L^2\cap L^{\infty}$ such that
\begin{align}
\lim_{t\to+\infty}\|\alpha_j(t)-\alpha_j^+\|_{L^2}=0
\label{L2convergence}
\end{align}
for $j=1,2$. Moreover we have 
$\alpha_1^+(\xi) \cdot \alpha_2^+(\xi)=0$ for $\xi \in \R$. 
\end{prp}

Once this proposition is obtained, we can derive Theorem~\ref{thm_initial} 
immediately by setting $\varphi_j^+=\op{F}^{-1}\alpha_j^+$. 
Indeed we have
\[
 \|u_j(t)-\op{U}(t)\varphi_j^+\|_{L^2}
=
 \|\op{F}(\op{U}(-t)u_j(t)-\varphi_j^+)\|_{L^2}
=
\|\alpha_j(t)-\alpha_j^+\|_{L^2}\to 0
\]
as $t\to +\infty$.

In the rest of this section, we will prove Proposition~\ref{prp_asym_prof}. 
Note that many parts of the arguments below are simliar to those in \cite{NS}, 
though we need several modifications to fit for the present situation. 
Before going into the proof, let us recall the following lemmas.

\begin{lem}
\label{lem_M}
Let $C_0>0$, $C_1\ge 0$, $p>1$, $q>1$ and $t_0\ge2$. 
Suppose that $\Phi(t)$ satisfies 
\begin{align*}
\frac{d\Phi}{dt}(t)
\le -\frac{C_0}{t}\left|\Phi(t)\right|^p + \frac{C_1}{t^{q}}
\end{align*}
for $t\ge t_0$. Then we have
\begin{align*}
\Phi(t)\le \frac{C_2}{(\log t)^{p^{*}-1}}
\end{align*}
for $t\ge t_0$, where $p^{*}$ is the H\"older conjugate of $p$ (i.e., 
$1/p+1/p^{*}=1$), and
\begin{align*}
C_2=\frac{1}{\log 2}\left( (\log t_0)^{p^{*}}\Phi(t_0) 
  + C_1\int_{2}^\infty\frac{(\log \tau)^{p^{*}}}{\tau^{q}}d\tau \right) 
  + \left( \frac{p^{*}}{C_0 p} \right)^{p^{*}-1}.
\end{align*}
\end{lem}
For the proof, see Lemma~4.1 in \cite{KMatsS}.

\begin{lem}\label{lem_ODE}
Let $t_0>0$ be given. For 
$\lambda$, $Q \in C\cap L^1([t_0,\infty))$, 
assume that $y(t)$ satisfies
\begin{align*}
\frac{dy}{dt}(t)
= \lambda(t)y(t) + Q(t)
\end{align*}
for $t\ge t_0$. Then we have 
\[
 |y(t)-y^+|
 \le 
 C_3 \int_t^\infty (|y^+||\lambda(\tau)|+|Q(\tau)|) \, d\tau
\]
for $t\ge t_0$, where 
\[
 C_3=\exp\left(\int_{t_0}^\infty |\lambda(\tau)|\, d\tau \right)
\]
and
\[
 y^+=
 y(t_0) e^{\int_{t_0}^\infty\lambda(\tau)\, d\tau}
 + 
 \int_{t_0}^\infty Q(s) e^{\int_{s}^\infty \lambda(\tau)\, d\tau}\, ds.
\]
\end{lem}
For the proof, see Lemma~4.2 in \cite{NS}.\\

\noindent{\em Proof of Proposition~\ref{prp_asym_prof}.}\ 
We first show the pointwise convergence of $\alpha(t,\xi)$ 
as $t\to +\infty$. 
We fix $\xi \in \R$ and introduce 
\[
 \rho(t,\xi)
 =
 2\realpart\Bigl[ 
  \cc{\alpha_1(t,\xi)}R_1(t,\xi)- \cc{\alpha_2(t,\xi)}R_2(t,\xi) 
  \Bigr].
\]
Then it follows from \eqref{profile_1} and \eqref{profile_2} that
\begin{align*}
 \pa_ t \Bigl(|\alpha_1(t,\xi)|^2- |\alpha_2(t,\xi)|^2\Bigr)
 =
 2\realpart\Bigl[ 
  \cc{\alpha_1}\pa_t \alpha_1- \cc{\alpha_2}\pa_t \alpha_2
  \Bigr]
 = \rho(t,\xi).
\end{align*}
Also \eqref{apriori2} and \eqref{R} lead to 
\begin{align*}
 \int _{2}^\infty\left| \rho(\tau,\xi)\right|\, d\tau 
 \le 
 C\int _{2}^\infty |\alpha(t,\xi)| |R(t,\xi)|\, d\tau 
 \le 
 \int _{2}^\infty 
 C \eps^4 \jb{\xi}^{-2}\tau^{3\gamma-5/4}\, 
 d\tau 
 \le
 C\eps^4 \jb{\xi}^{-2}
\end{align*}
for $0<\gamma< 1/12$.
Therefore we obtain 
\begin{align}\label{CL3}
|\alpha_1(t,\xi)|^2- |\alpha_2(t,\xi)|^2 
=
|\alpha_1(2,\xi)|^2- |\alpha_2(2,\xi)|^2
+
 \int _{2}^{t} \rho(\tau,\xi)\, d\tau 
=
m(\xi)-r(t,\xi),
\end{align}
where
\begin{align*}
m(\xi)
=
|\alpha_1(2,\xi)|^2- |\alpha_2(2,\xi)|^2
+
 \int _{2}^\infty \rho(\tau,\xi)\, d\tau 
\end{align*}
and
\begin{align*}
 r(t,\xi)=\int _{t}^\infty \rho(\tau,\xi)\, d\tau
\end{align*}
for $t \geq 2$.
Note that
\begin{align*}
 |m(\xi)|
 \le 
 |\alpha(2,\xi)|^2
 +\int _{2}^\infty|\rho(\tau,\xi)|\, d\tau 
\le 
 C\eps^2 \jb{\xi}^{-2} 
\end{align*}
and
\begin{align}
 |r(t,\xi)|\le 
\int _{t}^\infty |\rho(\tau,\xi)|\, d\tau
\le C\eps^4 \jb{\xi}^{-2}t^{3\gamma - 1/4}
\label{est_r}
\end{align}
for $0<\gamma< 1/12$.
Now we divide the argument into three cases according to the sign of 
$m(\xi)$ as follows.

\begin{itemize}
\item
\underline{\bf Case 1: $m(\xi)>0$.} \ 
First we focus on the asymptotics  for $\alpha_2$. 
By \eqref{CL3}, we can rewrite \eqref{profile_2} as
\[
\pa_t \alpha_2(t,\xi)=
-\frac{1}{t} |\alpha_2(t,\xi)|^2 \alpha_2(t,\xi) 
-\frac{m(\xi)}{t}\alpha_2(t,\xi)  +\frac{r(t,\xi)}{t}\alpha_2(t,\xi)
 +R_2(t,\xi)
\] 
for $t \geq 2$.
So we have 
\begin{align*}
 \pa_t\bigl( |\alpha_2(t,\xi)|^2\bigr)
=
 2\realpart \Bigl(\cc{\alpha_2}\pa_t \alpha_2\Bigr)
\le
 0 -\frac{2m(\xi)}{t} |\alpha_2(t,\xi)|^2 
 + C\eps^4 \jb{\xi}^{-2}t^{3\gamma-5/4}
\end{align*}
for $t \geq 2$,
whence
\begin{align*}
\pa_t \left( t^{2m(\xi)} |\alpha_2(t,\xi)|^2 \right)
\le 
 C\eps^4 \jb{\xi}^{-2}t^{3\gamma + 2m(\xi)- 5/4}.
\end{align*}
Integration in $t$ leads to
\begin{align*}
t^{2m(\xi)} |\alpha_2(t,\xi)|^2 -2^{2m(\xi)} |\alpha_2(2,\xi)|^2
\le 
C\eps^4 \jb{\xi}^{-2}\int^t_{2}\tau^{3\gamma + 2m(\xi)- 5/4}d\tau 
\le C\eps^4 \jb{\xi}^{-2}
\end{align*}
for $t\ge 2$. Therefore we see that 
\begin{align}
 |\alpha_2(t,\xi)| \le C\eps \jb{\xi}^{-1} t^{-m(\xi)}.
\label{decay_alpha2}
\end{align}
In particular, $\alpha_2(t,\xi) \to 0$ as $t\to +\infty$. 
Next we turn our attentions to the asymptotics for $\alpha_1$. 
Since \eqref{profile_1} can be viewed as 
\[
\pa_t \alpha_1(t)=\lambda(t) \alpha_1(t) +Q(t)
\] 
with $\lambda(t)=-|\alpha_2(t,\xi)|^2/{t}$ and $Q(t)=R_1(t,\xi)$,
we can apply Lemma~\ref{lem_ODE} to obtain
\[
 |\alpha_1(t,\xi)-\beta_1^+(\xi)|
 \le
 C\int_t^\infty \left(\frac{|\beta_1^+(\xi)| |\alpha_2(\tau,\xi)|^2}{\tau}
 + 
 |R_1(\tau,\xi)|\right)\, d\tau
\]
for $t \geq 2$,
where
\begin{align*}
\beta_1^{+}(\xi)
=
\alpha_1(2,\xi)
e^{-\int_{2}^{\infty} |\alpha_2(\tau,\xi)|^2\frac{d\tau}{\tau}}
+
\int_{2}^\infty 
R_1(s,\xi)
e^{-\int_s^\infty |\alpha_2(\tau,\xi)|^2\frac{d\tau}{\tau}}\, ds.
\end{align*}
By \eqref{apriori2}, \eqref{R} and \eqref{decay_alpha2}, we have
\begin{align}\label{est_beta1}
|\beta_1^{+}(\xi)|
&\le 
|\alpha_1(2,\xi)| 
+ \int_{2}^\infty |R_1(s,\xi)| ds 
\le 
 C\eps \jb{\xi}^{-1}
\end{align}
and
\begin{align*}
\int_t^\infty \left(\frac{|\beta_1^+(\xi)| |\alpha_2(\tau,\xi)|^2}{\tau}
 + 
 |R_1(\tau,\xi)|\right)\, d\tau
&\le 
C\int_t^\infty
\left(
 \frac{\eps^3 \jb{\xi}^{-3}}{\tau^{1+2m(\xi)}}
 +
 \frac{\eps^3 \jb{\xi}^{-1}}{\tau^{5/4-3\gamma}}
\right)\, d\tau\\
&\le
\frac{C\eps^3 \jb{\xi}^{-3}}{2m(\xi)t^{2m(\xi)}}
+
\frac{C\eps^3 \jb{\xi}^{-1}}{t^{1/4-3\gamma}}
\end{align*}
for $t \geq 2$.
Therefore we conclude that $\alpha_1(t,\xi) \to \beta_1^{+}(\xi)$ 
as $t\to +\infty$.

\item
\underline{\bf Case 2: $m(\xi)<0$.}\ 
Similarly to the previous case, we have
\begin{align*}
\lim _{t\to +\infty}|\alpha_1(t,\xi)|=0,
\qquad 
\lim _{t\to +\infty}|\alpha_2(t,\xi)-\beta_2^{+}(\xi)|=0 
\end{align*}
for each fixed $\xi \in \R$, where 
\begin{align*}
\beta_2^{+}(\xi)
:=
\alpha_2(2,\xi)
e^{-\int_{2}^\infty |\alpha_1(\tau,\xi)|^2\frac{d\tau}{\tau}}
+
\int_{2}^\infty 
R_2(s,\xi)
e^{-\int_{s}^\infty |\alpha_1(\tau,\xi)|^2\frac{d\tau}{\tau}}\, ds.
\end{align*}
Remark that $|\beta_2^+(\xi)|\le C\eps \jb{\xi}^{-1}$.

\item
\underline{\bf Case 3: $m(\xi)=0$.}\ 
By \eqref{profile_1}, \eqref{apriori2}, \eqref{R}, 
\eqref{CL3} and \eqref{est_r}, we have
\begin{align*}
\pa_t \left( |\alpha_1(t,\xi)|^2\right)
&
 \leq -\frac{2}{t} |\alpha_1(t,\xi)|^4 
- \frac{2r(t,\xi)}{t} |\alpha_1(t,\xi)|^2 
 + 2|\alpha_1(t,\xi)||R_1(t,\xi)| \\
&\le 
- \frac{2}{t} |\alpha_1(t,\xi)|^4 
 + C\eps^4 \jb{\xi}^{-2} t^{3\gamma - 5/4}
\end{align*}
for $t \ge 2$, and $0<\gamma< 1/12$. 
Thus we can apply Lemma~\ref{lem_M} with $\Phi(t)=|\alpha_1(t,\xi)|^2$ 
to obtain
\begin{align*}
 |\alpha_1(t,\xi)|
\le \frac{C}{(\log t)^{1/2}}\to 0 \qquad (t\to +\infty).
\end{align*}
Also \eqref{CL3} gives us 
$|\alpha_2(t,\xi)|=\sqrt{|\alpha_1(t,\xi)|^2 +r(t,\xi)} \to 0$ 
as $t\to +\infty$.
\end{itemize}

Summing up the three cases above, we deduce that 
$\alpha(t,\xi)$ converges as $t\to +\infty$ for each fixed $\xi\in \R$.
To obtain \eqref{L2convergence}, we set
\begin{align*}
\alpha^{+}_1(\xi):=
\left\{\begin{array}{cl}
  \beta^{+}_1(\xi) & \left(m(\xi)>0 \right),\\
  0 & \left(m(\xi)\le 0 \right),
\end{array}\right.
\quad
\alpha^{+}_2(\xi):=
\left\{ \begin{array}{cl}
 0 & \left( m(\xi) \ge0  \right),\\
 \beta^{+}_2(\xi) & \left( m(\xi) < 0 \right),
\end{array}\right.
\end{align*}
and $\alpha^{+}(\xi)=(\alpha^{+}_1(\xi), \alpha^{+}_2(\xi))$ 
for $\xi\in \R$, where $\beta_{1}^{+}(\xi)$ and $\beta_{2}^{+}(\xi)$
are shown in Cases 1 and 2, respectively.
Then it is obvious that 
$\alpha_1^+(\xi) \cdot \alpha_2^+(\xi)=0$ for $\xi \in \R$. 
Also, by virtue of \eqref{est_beta1}, we have 
$\alpha^{+}\in L^2\cap L^{\infty}(\R)$ 
and
\begin{align*}
\left| \alpha(t,\xi)-\alpha^{+}(\xi) \right|^2
\le 
C\eps^2 \jb{\xi}^{-2}\in L^1(\R)
\end{align*}
for $t\ge 2$. Moreover, it holds that
\[
 \lim_{t\to +\infty}
 \left| \alpha(t,\xi)-\alpha^{+}(\xi) \right|^2
 =0
\]
for each fixed $\xi \in \R$. 
Therefore  Lebesgue's dominated convergence theorem yields 
\eqref{L2convergence}.
\qed

\appendix\def\thesection{A}
\section{Proof of Lemma~\ref{remainder}}  
\label{sec_app_a}
We give a proof of Lemma~\ref{remainder}. 
For this purpose, we introduce some notations. We define the operators 
$\op{M}(t)$, $\op{D}(t)$ and $\op{W}(t)$ by
\[
 \bigl( \op{M}(t) \phi \bigr)(x)
 =
 e^{i\frac{x^2}{2t}}\phi(x),
 \quad
 \bigl( \op{D}(t) \phi \bigr)(x) 
 = 
 (it)^{-1/2} \phi \left(\frac{x}{t} \right),
 \quad 
  \op{W}(t)  \phi 
 =
\op{F} \op{M}(t) \op{F}^{-1} \phi,
\]
so that $\op{U}(t)=\exp(i\frac{t}{2}\pa_x^2)$ is decomposed into 
\begin{align}\label{MDFM}
\op{U}(t)=\op{M}(t) \op{D}(t) \op{F} \op{M}(t)
=\op{M}(t) \op{D}(t) \op{W}(t) \op{F}. 
\end{align}
An important estimate is
\begin{align}
 \|(\op{W}(t)  -1)\phi\|_{L^{\infty}} 
+
 \|(\op{W}(t)^{-1} -1)\phi\|_{L^{\infty}} 
 \le 
 Ct^{-1/4}\|\phi\|_{H^1},
\label{est_W}
\end{align}
which comes from the Gagliardo-Nirenberg inequality $\|\phi\|_{L^{\infty}}\le C\|\phi\|_{L^2}^{1/2}\|\pa_x \phi\|_{L^2}^{1/2}$ and the inequality
\begin{align}\label{ineq_key1}
|e^{i\theta} -1|\le C|\theta|^{\sigma} 
\qquad (\theta\in\R, \ 0\le\sigma\le1)
\end{align}
with $\theta=x^2/(2t),\sigma=1/2$. 
Note also that 
\begin{align}\label{ineq_key2}
 \|\op{W}(t)\op{F}\op{U}(-t)\phi\|_{H^1} 
+ \|\op{W}(t)^{-1}\op{F}\op{U}(-t)\phi\|_{H^1} 
 \le C(\|\phi\|_{L^2}+\|\op{J}\phi\|_{L^2})
\end{align}
and
\begin{align}\label{ineq_key3}
 \|\op{F} \op{U}(-t) [ \phi_1 \phi_2 \phi_3 ] \|_{L^{\infty}}
 \le
 C\|\phi_1\|_{L^2} \|\phi_2\|_{L^2} \|\phi_3\|_{L^{\infty}},
\end{align}
where the constant $C$ is independent of $t$ 
(see e.g., \cite{LiS1} for the proof). 
In what follows, we will occasionally omit ``$(t)$" 
from $\op{M}(t)$, $\op{D}(t)$, $\op{W}(t)$ 
if it causes no confusion. 

Let $\alpha$ be given by \eqref{def_alpha}. 
By \eqref{MDFM}, we have
\begin{align*}
\op{F}\op{U}(-t)\bigl[ |u_2|^2u_1\bigr]
=
\op{W}^{-1}\op{D}^{-1}\op{M}^{-1}
\bigl[ |\op{M}\op{D}\op{W}\alpha_2|^2 \op{M}\op{D}\op{W}\alpha_1\bigr]
=
\frac{1}{t}\op{W}^{-1}\bigl[ |\op{W}\alpha_2|^2 \op{W}\alpha_1\bigr],
\end{align*}
whence
\begin{align*}
 R_1
&=
\frac{1}{t}
\Bigl( 
|\alpha_2|^2 \alpha_1
-
\op{W}^{-1}\bigl[ |\op{W}\alpha_2|^2 \op{W}\alpha_1\bigr]
\Bigr)\\
&=
\frac{1}{t}
(1-\op{W}^{-1})\bigl[ |\op{W}\alpha_2|^2 \op{W}\alpha_1\bigr]
+
\frac{1}{t}|\op{W}\alpha_2|^2 (1-\op{W})\alpha_1\\
&\quad+
\frac{1}{t}(\op{W}\alpha_2)(\cc{(1-\op{W})\alpha_2}) \alpha_1
+\frac{1}{t}((1-\op{W})\alpha_2) \cc{\alpha_2} \alpha_1.
\end{align*}
Therefore \eqref{est_W}, \eqref{ineq_key2}, \eqref{ineq_key3}  and 
the Sobolev imbedding $H^1(\R) \hookrightarrow L^{\infty}(\R)$ lead to 
\begin{align}
 |R_1(t,\xi)|
 \le
 Ct^{-5/4} (\|u\|_{L^2}+\|\op{J}u\|_{L^2})^3.
\label{est_R_0}
\end{align}
Next we observe that
\begin{align*}
i\xi R_1
&=
\frac{i\xi}{t}
|\alpha_2|^2 \alpha_1
-
\op{F}\op{U}(-t)\Bigl[ \pa_x \bigl(|u_2|^2 u_1\bigr)\Bigr]
\\
&=
\frac{1}{t}
\Bigl( 
\alpha_2^{(1)}\cc{\alpha_2} \alpha_1
-
\op{W}^{-1}
\bigl[ (\op{W}\alpha_2^{(1)})(\cc{\op{W}\alpha_2}) \op{W}\alpha_1\bigr]
\Bigr)
\\
&\quad+
\frac{1}{t}
\Bigl(
\alpha_2 \cc{\alpha_2^{(1)}} \alpha_1
-
\op{W}^{-1}
\bigl[ (\op{W}\alpha_2)(\cc{\op{W}\alpha_2^{(1)}}) \op{W}\alpha_1\bigr]
\Bigr)
\\
 &\quad+
\frac{1}{t}
\Bigl( 
\alpha_2\cc{\alpha_2} \alpha_1^{(1)}
-
\op{W}^{-1}
\bigl[ (\op{W}\alpha_2)(\cc{\op{W}\alpha_2}) \op{W}\alpha_1^{(1)}\bigr]
\Bigr),
\end{align*}
where $\alpha_j^{(1)}=i\xi \alpha_j$. Then we see as before that  
\begin{align}
|\xi R_1(t,\xi)| \le
  Ct^{-5/4} (\|u\|_{H^1}+\|\op{J}u\|_{H^1})^3.
\label{est_R_1}
\end{align}
The desired estimate for $R_1$ follows immediately from \eqref{est_R_0} and 
\eqref{est_R_1}. The estimate for $R_2$ can be shown in the same way.
\qed

\appendix\def\thesection{B}
\section{Final state problem  for \eqref{eq}}  
\label{sec_app_b}

To complement Remark~\ref{rem1}, we give two auxiliary results on the final 
state problem for \eqref{eq}, that is, finding a solution $u=(u_1,u_2)$ to 
\eqref{eq} which satisfies
\begin{align}
 \lim_{t\to +\infty} \|u_j(t)-\op{U}(t)\psi_j^+\|_{L^2}=0,\quad j=1,2,
\label{scattering}
\end{align}
for a prescribed final state $\psi^+=(\psi_1^+,\psi_2^+)$. Roughly speaking, 
the propositions below imply that 
\eqref{scattering} holds if and only if 
\begin{align}\label{decoupling_fin}
\hat{\psi}_1^+(\xi)\cdot \hat{\psi}_2^+(\xi)=0,\quad \xi \in \R.
\end{align}
This indicates that our problem must be distinguished from the usual 
short-range situation, because \eqref{scattering} should hold 
in the short-range case regardless of 
whether \eqref{decoupling_fin} is true or not (see e.g. \cite{HLO}). 

The precise statements are as follows.

\begin{prp}\label{prp_final_a}
Let $T_0\geq 1$ be given, and let $u$ be a solution to \eqref{eq} for 
$t\ge T_0$ satisfying 
\begin{align}\label{assump_u}
\sup_{t\ge T_0}
\Bigl(
 t^{-\gamma} \|\op{U}(-t)u(t)\|_{H^{1,1} }
 + 
 \|\op{F}\op{U}(-t)u(t)\|_{L^{\infty}}
\Bigr)
<\infty
\end{align}
with some $\gamma\in(0, 1/12)$.
If there exists $\psi^+\in L^2$ with 
$\hat{\psi}^+\in L^{\infty}$ such that \eqref{scattering} holds, 
then we must have \eqref{decoupling_fin}.
\end{prp}

\begin{prp}\label{prp_final_b}
Suppose that  $\psi^+ $ satisfies $\hat{\psi}^+\in H^{0,s}\cap L^{\infty}$ 
with some $s>1$, and that $\delta=\|\hat{\psi}^+\|_{L^{\infty}}$ is suitably small. 
If \eqref{decoupling_fin} holds, then there exist $T\ge 1$ and 
a unique solution $u$ to \eqref{eq} for $t\ge T$ 
satisfying $\op{U}(-t)u\in C([T,\infty);H^{0,1})$
and  \eqref{scattering}. 
\end{prp}

We are going to give a proof of them. 
Note that the arguments below are essentially the same as those given in 
Section~5 of \cite{Li}. We also remark that 
the regularity assumptions in these propositions are certainly not optimal. 
It may be possible to relax them (see e.g. \cite{HN2}), 
but that is out of the present purpose. \\

\noindent {\em Proof of Proposition~\ref{prp_final_a}.}\ 
In what follows, we write 
$N_1(v)=|v_2|^2v_1$, $N_2(v)=|v_1|^2v_2$ and 
$N(v)=(N_1(v),N_2(v))$ for $v=(v_1,v_2)$. 
Let $\alpha$ be given by \eqref{def_alpha} for the solution $u$ 
to \eqref{eq}.  
Then, similarly to \eqref{profile_1}, we have
\begin{align*}
\pa_t \alpha_j(t,\xi)
=
-\frac{1}{t} N_j(\hat{\psi}^+(\xi)) + S_j(t,\xi) + R_j(t,\xi), 
\quad j=1,2,
\end{align*}
where
\[
 S_j(t,\xi)
 = 
 \frac{1}{t} \left(N_j(\hat{\psi}^+(\xi)) - N_j(\alpha(t,\xi))  \right)
\]
and 
\[R_j(t,\xi)
=
\frac{1}{t} N_j(\alpha(t,\xi))
-\op{F}\Bigl[\op{U}(-t) N_j(u(t,\cdot)) \Bigr](\xi).
\] 
Now we shall argue by contradiction. If \eqref{decoupling_fin} is not true, 
then we can take $\eta>0$ such that 
$\|N_j(\hat{\psi}^+)\|_{L^2}\ge \eta$ for $j=1,2$.
By \eqref{assump_u} and Lemma~\ref{remainder}, we have 
$\|R_j(t)\|_{L^2} \le Ct^{-5/4+3\gamma}$ for $t\ge T_0$. 
We also note that 
\begin{align*}
 \|S_j(t)\|_{L^2} 
 &=\frac{1}{t} \|N_j(\hat{\psi}^+) - N_j(\op{F}\op{U}(-t) u(t))\|_{L^2}\\
 &\le 
\frac{C}{t}
 \bigl(
  \|\hat{\psi}\|_{L^{\infty}}^2 +  \|\op{F}\op{U}(-t) u(t)\|_{L^{\infty}}^2
 \bigr)
 \|\op{F}(\psi^+- \op{U}(-t) u(t))\|_{L^2}\\
&\le
 \frac{C}{t} \|\op{U}(t) \psi^+ -u(t) \|_{L^2},
\end{align*}
whence, by \eqref{scattering}, we can take $T^* \ge T_0$ such that 
$\|S_j(t)\|_{L^2}\le \eta/(2t)$ for $t\ge T^*$. 
Summing up, we obtain
\begin{align*}
 \|\op{U}(-2t)u_j(2t)-\op{U}(-t)u_j(t)\|_{L^2}
 &= \|\alpha_j(2t)-\alpha_j(t)\|_{L^2} \\
 &\ge
 \eta \int_t^{2t}\frac{d\tau}{\tau} 
 -\int_t^{2t}\, \|S_j(\tau)\|_{L^2}d\tau
 -\int_t^{2t}\, \|R_j(\tau)\|_{L^2}d\tau\\
 &\ge
\frac{\eta}{2}\log 2 -Ct^{-1/4+3\gamma}
\end{align*}
for $t\ge T^{*}$. Letting $t\to \infty$, we have 
\[
 0=\|\psi_j^+ - \psi_j^+\|_{L^2}\ge \frac{\eta}{2}\log 2>0,
\]
which is the desired contradiction. 
\qed\\

\noindent {\em Proof of Proposition~\ref{prp_final_b}.} \ 
With $T\ge 1$ to be fixed, let us introduce the function space 
\[
 \mathfrak{X}_T=
\Bigl\{
 \phi=(\phi_1(t,x),\phi_2(t,x))
 \, \Bigm|\, 
 \op{U}(-t) \phi(t,\cdot) \in C([T, \infty);H^{0,1}) 
\Bigr\}
\]
and the norm
\[
 \|\phi\|_{\mathfrak{X}_T}
 =
 \sup_{t\in [T, \infty)}\bigl(
  t^{\mu+1/2} \|\phi(t)\|_{L^2} +t^{\mu} \|\op{J} \phi(t)\|_{L^2}
 \bigr),
\]
where $0<\mu<(s_{0}-1)/2$ and $s_{0}=\min\{2,s\}$. For $v=(v_1,v_2)\in \mathfrak{X}_T$, 
we set 
\begin{align}
 \Phi_j[v](t)
 =
 \op{U}(t) \psi_j^+
  -\int_t^\infty \op{U}(t-\tau)N_j(v(\tau))\, d\tau,
 \quad j=1,2,
 \label{def_Phi}
\end{align}
and $\Phi[v]=(\Phi_1[v], \Phi_2[v])$. 
We also put 
$w^\sharp (t)=\op{M}(t)\op{D}(t) \op{F}\psi^+$,  
$w^{\flat}(t)=\op{U}(t)\psi^+ -w^\sharp(t)$, 
$\kappa=\|\psi^+\|_{H^{0,s_0}}$
and 
\[
 \mathfrak{Y}_{T}
= \Bigl\{ 
 \phi \in \mathfrak{X}_T\, \Bigm| \, \|\phi-w^\sharp\|_{\mathfrak{X}_T}\le \kappa
\Bigr\}.
\]
Since \eqref{decoupling_fin} yields $N(\hat{\psi}^+) =0$, it follows from 
\eqref{MDFM} that 
\begin{align}
 \op{U}(-\tau) N(w^\sharp(\tau))
 =
 \frac{1}{\tau}\op{M}(\tau)^{-1}\op{F}^{-1} N(\hat{\psi}^+) 
 =0.
\label{vanish}
\end{align}
We observe the basic estimates for $w^\sharp (t)$ and $w^\flat (t)$: 
\begin{align*}
&\|w^\sharp (t)\|_{L^{\infty}}
 = t^{-1/2}\|\hat{\psi}^+\|_{L^{\infty}}
 = \delta t^{-1/2},\\
&\|w^\sharp (t)\|_{L^{2}}
 =\|\psi^+\|_{L^2}
 \le \kappa,\\
&\|w^\flat (t)\|_{L^2}
 \le \|(\op{M}(t)-1)\psi^+\|_{L^2}
 \le Ct^{-s_{0}/2}\|\psi^+\|_{H^{0,s_0}}
 \le C\kappa t^{-s_{0}/2},\\
&\|\op{J} w^\sharp (t) \|_{L^{2}}
 = \|x\psi^+\|_{L^2}
 \le \kappa,\\
&\|\op{J} w^\flat (t)\|_{L^2}
 \le \|x(\op{M}(t)-1)\psi^+\|_{L^2}
 \le Ct^{-(s_{0}-1)/2}\|\psi^+\|_{H^{0,s_0}}
 \le C\kappa t^{-(s_{0}-1)/2},
\end{align*}
where we have used  \eqref{J_to_U}, \eqref{MDFM} and 
\eqref{ineq_key1} with $\sigma=s_{0}/2$ or $(s_{0}-1)/2$.

Now we are going to show that $\Phi$ is a contraction mapping on 
$\mathfrak{Y}_{T}$ by choosing $\delta$ and $T$ appropriately. 
Let $v\in \mathfrak{Y}_{T}$. By using \eqref{vanish}, 
we rewrite \eqref{def_Phi} as 
\[
 \Phi[v](t) - w^\sharp (t)
 = 
 -\int_t^\infty\op{U}(t-\tau)
 \Bigl( N(v(\tau)) -N(w^\sharp (\tau)) \Bigr)\, d\tau
 +w^\flat(t).
\]
It follows from \eqref{kla_ineq} that 
\[
\|v(t) -w^\sharp (t) \|_{L^{\infty}}
 \le
 Ct^{-1/2} 
 \|v(t) -w^\sharp (t) \|_{L^2}^{1/2} 
 \|\op{J}(v(t) -w^\sharp (t) )\|_{L^2}^{1/2} 
 \le C\kappa t^{-3/4 -\mu}.
\]
So we have 
\begin{align*}
 \|v(t) \|_{L^{\infty}} 
 \le 
  \|w^\sharp (t) \|_{L^{\infty}}
 +
 \|v(t) -w^\sharp (t) \|_{L^{\infty}}
 \le
 C (\delta+  \kappa T^{-1/4-\mu})t^{-1/2}.
\end{align*}
Therefore 
\begin{align}
 \|\Phi[v](t) -w^\sharp (t)\|_{L^2}
 &\le
 C\int_t^\infty
 (\|v (\tau) \|_{L^{\infty}}^2+\|w^\sharp (\tau) \|_{L^{\infty}}^2)
 \|v(\tau) -w^\sharp (\tau) \|_{L^2}\, d\tau
 + 
 C\kappa t^{-s_{0}/2}
 \notag\\
 &\le
 C(\delta+  \kappa T^{-1/4-\mu})^2\kappa 
 \int_t^\infty \frac{d\tau}{\tau^{3/2+\mu}}
 + 
 C\kappa T^{-(s_{0}-1)/2+\mu} t^{-1/2-\mu}
 \notag\\
 &\le
 C(\delta^2+\kappa^2 T^{-1/2-2\mu} + T^{-(s_{0}-1)/2+\mu})\kappa t^{-1/2-\mu}.
 \label{est_contra0}
\end{align}
Also, because of the estimate 
\begin{align*}
 \|\op{J}(N(v(t))-N(w^\sharp (t)))\|_{L^2}
 \le&
 C(\|v\|_{L^{\infty}}^2+\|w^\sharp \|_{L^{\infty}}^2)
 \|\op{J}(v-w^\sharp )\|_{L^2}
\\
 &+ 
 C(\|v\|_{L^{\infty}}+\|w^\sharp \|_{L^{\infty}})
 (\|\op{J}v\|_{L^2}+\|\op{J} w^\sharp \|_{L^2}) \|v-w^\sharp \|_{L^{\infty}}
\\
 \le&
 C(\delta^2 + \kappa^2 T^{-1/2-2\mu}) \kappa t^{-1-\mu}\\
 &+
 C(\delta+\kappa T^{-1/4-\mu})t^{-1/2}\cdot 
 C(\kappa t^{-\mu}+\kappa ) \cdot  C\kappa t^{-3/4-\mu}\\
 \le&
 C(\delta^2 + \kappa^2 T^{-1/2-\mu}+\delta\kappa T^{-1/4} )\kappa t^{-1-\mu},
\end{align*}
we obtain
\begin{align}
 \|\op{J} (\Phi[v](t) -w^\sharp (t))\|_{L^2}
 &\le
 \int_t^\infty 
 \|\op{J}(N(v(\tau))-N(w^\sharp (\tau)))\|_{L^2}\, d\tau
 + 
 C\kappa t^{-(s_{0}-1)/2}
 \notag\\
 &\le
 C(\delta^2 + \kappa^2 T^{-1/2-\mu}+\delta\kappa T^{-1/4} )\kappa 
 \int_t^\infty \frac{d\tau}{\tau^{1+\mu}}
 + 
 C\kappa T^{-(s_{0}-1)/2+\mu} t^{-\mu}
 \notag\\
 &=
 C(\delta^2 + \kappa^2 T^{-1/2-\mu}+\delta\kappa T^{-1/4} + T^{-(s_{0}-1)/2+\mu}) 
\kappa t^{-\mu}.
 \label{est_contra1}
\end{align}
Combining \eqref{est_contra0} and \eqref{est_contra1}, we arrive at 
\begin{align*}
 \|\Phi[v] -w^\sharp \|_{\mathfrak{X}_T}
 \le
 \underbrace{
 C(\delta^2 + \kappa^2 T^{-1/2-\mu}+\delta\kappa T^{-1/4}+ T^{-(s_{0}-1)/2+\mu})
}_{(*)} \kappa.
\end{align*}
Hence we have $\Phi[v] \in \mathfrak{Y}_{T}$ 
if we choose $\delta$ so small and $T$ so large 
that the term $(*)$ does not exceed $1$.
Next we take $v$, $\tilde {v} \in \mathfrak{Y}_{T}$. 
Then we have 
\[
\Phi[v](t) -\Phi[\tilde{v}](t)
=
-\int_t^\infty
 \op{U}(t-\tau)\Bigl(N(v(\tau))-N(\tilde{v}(\tau))\Bigr)
\, d\tau
\]
and we can see as before that 
\[
 \|\Phi[v]-\Phi[\tilde{v}]\|_{\mathfrak{X}_T} 
 \le 
 \frac{1}{2} \|v-\tilde{v}\|_{\mathfrak{X}_T}
\]
by choosing $\delta$ and $T$ suitably. 
Therefore $\Phi:\mathfrak{Y}_{T}\to \mathfrak{Y}_{T}$ 
is a contraction mapping, and thus, admits a unique fixed point. 
In other words, there exists 
$u\in \mathfrak{Y}_{T}$ such that 
\[
 u(t)=\op{U}(t)\psi^+ -\int_t^\infty \op{U}(t-\tau)N(u(\tau))\, d\tau, 
\]
which gives the desired solution to \eqref{eq} for $t\ge T$. 
Moreover we have
\begin{align*}
 \|u(t)-\op{U}(t)\psi^+\|_{L^2}
 &\le
 \|u(t)-w^\sharp (t)\|_{L^2}
 + 
 \|w^{\flat}(t)\|_{L^2}\\
 &\le
 \kappa t^{-1/2-\mu} +C\kappa t^{-s_{0}/2}\\
 &\to  0
\end{align*}
as $t\to +\infty$. This completes the proof of Proposition \ref{prp_final_b}.
\qed

\medskip
\subsection*{Acknowledgments}
The work of H.~S. is supported by Grant-in-Aid for Scientific Research (C) 
(No.~17K05322), JSPS.


\end{document}